\documentclass{article}
\usepackage{amsmath,amsfonts,amssymb}
\usepackage[latin2]{inputenc}                                                   
\newtheorem{theorem}{Theorem}[section]
\newtheorem{definition}[theorem]{Definition}

\newtheorem{lemma}[theorem]{Lemma}

\newtheorem{remark}[theorem]{Remark}

\newcommand{\mC}{\mathbb C}

\newcommand{\mN}{\mathbb N}

\newcommand{\mR}{\mathbb R}

\newcommand{\be}{\begin{eqnarray}}
\newcommand{\ee}{\end{eqnarray}}
\newcommand{\bd}{\begin{definition}}
\newcommand{\ed}{\end{definition}}

\newcommand{\br}{\begin{remark}}
\newcommand{\er}{\end{remark}}
\newcommand{\gog}{{\mathfrak g}}
\newcommand{\bt}{\begin{tabular}}
\newcommand{\et}{\end{tabular}}

\def\Sp{\mathop{\rm Sp}\nolimits}
\def\SL{\mathop{\rm SL}\nolimits}

\def\Mp{\mathop{\rm Mp}\nolimits}
\def\End{\mathop{\rm End}\nolimits}

\def\sp{\mathop{\mathfrak{sp}}\nolimits}
\def\sl{\mathop{\mathfrak{sl}}\nolimits}
\def\mp{\mathop{\mathfrak{mp}}\nolimits}

\def\Id{\mathop{\rm Id}\nolimits}

\input xypic
\input epsf
\input xy
\xyoption{all}

\begin{document}

\date{}

\baselineskip13pt

\title{Symplectic twistor operator and its solution space on the standard symplectic space $(\mR^2,\omega)$}
\author{Marie Dost\'alov\'a, Petr Somberg
}

\maketitle

\abstract 
We introduce the symplectic twistor operator $T_s$ in the symplectic spin geometry
of real dimension two, as a symplectic analogue of the Dolbeault operator in the
complex spin geometry of complex dimension $1$. Based on the techniques of the 
metaplectic Howe duality and the algebraic Weyl algebra, we compute the space of its 
solutions on the standard symplectic space $(\mR^2,\omega)$. 


\endabstract


\section{Introduction and Motivation}

Central problems and questions in differential geometry of Riemannian spin manifolds 
are usually reflected in analytic and spectral properties 
of the pair of first order differential operators acting on spinors, the Dirac operator and the twistor 
operator. In particular, there is rather subtle relation between geometry and
topology of a given manifold and the spectra resp. the solution spaces of these operators, 
see e.g., \cite{thom}, \cite{bfkg} and references therein.   

Based on the Segal-Shale-Weil representation, the symplectic version of the Dirac operator $D_s$ 
was introduced in \cite{KOS}, and some of its basic analytic and spectral properties were
studied in \cite{crum}, \cite{MR2252919}, \cite{KAD}. 
Introducing the metaplectic Howe duality, \cite{bss}, a representation theoretical characterization
of the solution space of the symplectic Dirac operator was determined on the standard
symplectic space $(\mR^{2n},\omega)$. However, an explicit analytic description of this space 
is still missing and this fact has also substantial consequences for the 
present article.    

A variant of the first order symplectic twistor operator $T_s$ was introduced in $\cite{KAD}$ in 
the framework of the contact parabolic geometry, inducing the symplectic twistor operator 
on the symplectic leaves of foliation. Several basic properties including the solution space of the 
symplectic twistor operator on $\mR^{2n}$ are discussed in $\cite{ds}$. In particular, the
case $n=1$ relates to the framework developed in $\cite{ds}$ as well, but all the results 
for $n=1$ and $n> 1$ are intrinsically different. We remark that the approach in
$\cite{ds}$ is based on the procedure of geometrical prolongation of 
the symplectic twistor differential equation. The problem behind the 
case $n=1$ is that many first order operators (e.g., the Dirac and twistor operators on spinors) 
coincide in the case of one complex dimension
with the Cauchy-Riemann (Dolbeault) and its conjugate operators.

The aim of the present article is to fill this gap and discuss the case of $n=1$ by different 
methods, namely, by analytical and combinatorial techniques. A part of
the problem of finding the solution space of $T_s$ is the discovery of certain canonical 
representative solutions of the symplectic Dirac operator $D_s$ and the discovery of
certain non-trivial identities in the algebraic Weyl algebra.

 The system of partial differential equations representing $T_s$ is overdetermined,
acting on the space of functions valued in an infinite dimensional vector space of the 
Segal-Shale-Weil 
representation, and the solution space of $T_s$ is (even locally) infinite dimensional.
 Notice that the techniques of the metaplectic Howe duality are not restricted to $(\mR^2,\omega)$, 
but it is not straightforward for $({\mR^{2n},\omega})$, $n>1$, to write more explicit formulas for 
solutions with values in the higher dimensional non-commutative algebraic Weyl algebra.

The structure of our article goes as follows. In the first Section, we review basic properties 
of the symplectic spin geometry in the real dimension $2$, with emphasis on the metaplectic Howe duality. In 
the second Section, we give a general definition of the symplectic twistor operator $T_s$. The space 
of polynomial solutions of $T_s$ on $(\mR^2,\omega)$ is analysed in Section three, relying on two basic principles.
The first one is representation theoretical, coming from the action of the metaplectic Lie algebra on the function 
space of interest. The second one is then the construction of representative solutions in the 
particular irreducible subspaces of the function space. As a byproduct of our approach, we construct 
specific polynomial solutions of the symplectic Dirac operator $D_s$, which is 
also a novelty. In the end of this section, we indicate the collection of unsolved problems directly 
related to the topic of the present article.  

Throughout the article, we use the notation $\mN_0$ for the set of natural numbers including zero
and $\mN$ for the set of natural numbers without zero.


\subsection{Metaplectic Lie algebra $\mp(2,\mR)$, symplectic Clifford algebra and a class of simple 
weight modules for $\mp(2,\mR)$}

In the present section we recollect some basic algebraic and representation theoretical information needed 
in the analysis of the solution space of the symplectic twistor operator $T_s$, see e.g., \cite{bss}, \cite{crum}, 
\cite{fh}, \cite{MR2252919}, \cite{KAD}.

Let us consider a $2$-dimensional symplectic vector space $(\mR^{2},\omega=dx\wedge dy )$, and a 
symplectic basis $\{e, f\}$ with respect to the non-degenerate two form $\omega\in\wedge^2(\mR^{2})^\star$.
The linear action of $\sp(2,\mR)\simeq \sl(2,\mR)$ on $\mR^2$ induces the action on its tensor representations,  
and we have $g^\star\omega=\omega$ for all $g\in \sp(2,\mR)$. The set of three matrices 
$$H=\begin{pmatrix}
1 & 0 \\
0 & -1 
\end{pmatrix},\,
X=\begin{pmatrix}
0 & 1 \\
0 & 0 
\end{pmatrix},\,
Y=\begin{pmatrix}
0 & 0 \\
1 & 0 
\end{pmatrix}$$
is a basis of $\sp(2,\mR)$. 

The metaplectic Lie algebra $\mp(2,\mR)$ is the Lie algebra of the twofold group covering 
$\pi : \Mp(2,\mR)\to \Sp(2,\mR)$ of the symplectic Lie group $\Sp(2,\mR)$. It can be realized 
by homogeneity two elements in 
the symplectic Clifford algebra $Cl_s(\mR^{2},\omega)$, where the homomorphism $\pi_\star : \mp(2,\mR)\to \sp(2,\mR)$ 
is given by
\begin{eqnarray}
& & \pi_\star (e\cdot e)=-2X,
\nonumber \\
& & \pi_\star (f\cdot f)=2Y,
\nonumber \\
& & \pi_\star (e\cdot f+f\cdot e)=2H.
\end{eqnarray}
\begin{definition}
The symplectic Clifford algebra $Cl_s(\mR^{2},\omega)$ is an associative unital algebra over $\mC$,
realized as a quotient of the tensor algebra $T(e,f)$ by a
two-sided ideal $I\,\subset T(e,f)$, generated by
$$
v_i\cdot v_j-v_j\cdot v_i=-i\omega (v_i,v_j)
$$ 
for all $v_i,v_j\in\mR^{2}$.
\end{definition}
The symplectic Clifford algebra $Cl_s(\mR^{2},\omega)$ is isomorphic to the Weyl algebra $W_2$ of complex 
valued algebraic 
differential operators on $\mR$, and the symplectic Lie algebra $\sp(2,\mR)$ can be realized as 
a subalgebra of $W_2$. In particular, the Weyl algebra is an associative algebra generated by 
$\{q,\partial_q\}$, the multiplication operator by $q$ and, 
differentiation $\partial_q$, and the symplectic Lie algebra $\sp(2,\mR)$ has a basis $\{-\frac{i}{2} q^2,
-\frac{i}{2} \frac{\partial^2}{\partial q^2}, q \frac{\partial}{\partial q}+\frac{1}{2}\}$.

The symplectic spinor representation is the irreducible Segal-Shale-Weil representation of $Cl_s(\mR^{2},\omega)$ on 
$L^2(\mR,e^{-\frac{q^2}{2}} dq_{\mR})$, the space of square integrable functions on 
$(\mR,d\mu =e^{-\frac{q^2}{2}} dq_{\mR})$ with $dq_{\mR}$ the Lebesgue measure. 
Its action, the symplectic Clifford multiplication $c_s$, preserves the subspace
of $C^\infty$(smooth)-vectors given by the Schwartz space $S(\mR)$ of rapidly decreasing complex valued functions on
$\mR$ as its dense subspace. The space $S(\mR)$ can be regarded as a smooth Frechet globalization
of the space of $\tilde{K}=\widetilde{\mathrm{U}}(1)$-finite vectors in the representation, 
where $\tilde{K}\subset \Mp(2,\mR)$ is the maximal compact subgroup given by the double cover of 
$K=\mathrm{U}(1)\subset \Sp(2,\mR)$. Though we shall work in the smooth globalization $S(\mR)$, our representative
vectors constructed in Section $3$ belong to the underlying Harish-Chandra module
of $\tilde{K}=\widetilde{\mathrm{U}}(1)$-finite vectors preserved by $c_s$. 

The function spaces associated to Segal-Shale-Weil representation are supported on $\mR\,\subset\mR^2$, a 
maximal isotropic subspace of $(\mR^2,\omega)$. 
In its restriction to $\mp(2,\mR)$, it decomposes into two 
unitary representations realized on the subspace of even resp. odd functions:
\begin{eqnarray}\label{sshrepr}
\varrho :\mp(2,\mR)\to \End(S(\mR)),
\end{eqnarray} 
where the basis vectors act by
\begin{eqnarray}
& & \varrho (e\cdot e)=iq^2,
\nonumber \\
& & \varrho (f\cdot f)=-i\partial_q^2,
\nonumber \\
& & \varrho (e\cdot f+f\cdot e)=q\partial_q+\partial_qq.
\end{eqnarray}
In this representation $Cl_s(\mR^{2},\omega)$ acts 
on $L^2(\mR,e^{-\frac{q^2}{2}} dq_{\mR})$ by 
continuous unbounded operators with the domain $S(\mR)$. The space of
$\tilde{K}=\widetilde{\mathrm{U}}(1)$-finite vectors has a basis $\{q^je^{-\frac{q^2}{2}}\}_{j=0}^\infty$,
its even $\mp(2,\mR)$-submodule $\{q^{2j}e^{-\frac{q^2}{2}}\}_{j=0}^\infty$ resp.
odd $\mp(2,\mR)$-submodule $\{q^{2j+1}e^{-\frac{q^2}{2}}\}_{j=0}^\infty$.
It is also an irreducible representation of 
$\mp(2,\mR)\ltimes h(2)$, the semidirect product of $\mp(2,\mR)$ and a $3$-dimensional 
Heisenberg Lie algebra spanned by $\{e,f,\Id\}$. In the article we denote the Segal-Shale-Weil 
representation by 
${\mathcal S}$ and we have ${\mathcal S}\simeq {\mathcal S}_+\oplus{\mathcal S}_-$ as $\mp(2,\mR)$-module.

Let us denote by $\mathrm{Pol}(\mR^2)$ the vector space of complex valued polynomials on $\mR^2$, 
and by $\mathrm{Pol}_l(\mR^2)$ the subspace of homogeneity $l$ polynomials. The complex vector space
$\mathrm{Pol}_l(\mR^2)$ is as an irreducible $\mp(2,\mR)$-module isomorphic to $\mathrm{S}^l(\mC^{2})$, 
the $l$-th symmetric power of the complexification of the fundamental vector representation $\mR^2$, $l\in\mN_0$.
 

\subsection{Segal-Shale-Weil
representation and the metaplectic Howe duality}

Let us recall a representation-theoretical result of \cite{BL}. 
Let $\lambda_1$ be the fundamental weight of the Lie algebra $\sp(2,\mR)$, and let $L(\lambda)$ 
denote the simple module over universal
enveloping algebra $\mathcal{U}(\mp(2,\mR))$ of $\mp(2,\mR)$ generated by the highest weight vector of the weight 
$\lambda$.
Then the Segal-Shale-Weil representation for $\mp(2,\mR)$ is the highest weight representation  
$L(-\frac{1}{2}\lambda_1)\oplus L(-\frac{3}{2}\lambda_1)$. The highest weight vector is the eigenvector 
of the generator of $1$-dimensional maximal commutative subalgebra of $\mp(2,\mR)$.

The decomposition of the space of polynomial functions on $\mR^2$ valued in the Segal-Shale-Weil 
representation corresponds to the tensor product of 
$L(-\frac{1}{2}\lambda_1)\oplus L(-\frac{3}{2}\lambda_1)$ with 
symmetric powers $\mathrm{S}^l(\mC^{2n})$, $l\in\mN_0$, of the fundamental vector 
representation $\mC^{2}$ of $\sp(2,\mR)$. Note that all summands in the decomposition are 
again irreducible representations of $\mp(2,\mR)$.

\begin{lemma}(\cite{BL})\label{decomposition}
Let $\,l\in\mN_0$.
\begin{enumerate}
\item  
We have for $L(-\frac{1}{2}\lambda_1)$ and any $l$:
\begin{eqnarray*}
L(-\frac{1}{2}\lambda_1)\otimes \mathrm{S}^l(\mC^{2}) &\simeq &  
L(-\frac{1}{2}\lambda_1)\oplus L(\lambda_1-\frac{1}{2}\lambda_1)\oplus\dots 
\nonumber \\
& & \oplus L((l-1)\lambda_1-\frac{1}{2}\lambda_1)\oplus
L(l\lambda_1-\frac{1}{2}\lambda_1),
\end{eqnarray*}
\item 
We have for $L(-\frac{3}{2}\lambda_1)$ and any $l$:
\begin{eqnarray*}
L(-\frac{3}{2}\lambda_1)\otimes \mathrm{S}^l(\mC^{2}) &\simeq & 
L(-\frac{3}{2}\lambda_1)\oplus L(\lambda_1-\frac{3}{2}\lambda_1)\oplus\dots 
\nonumber \\
& & \oplus L((l-1)\lambda_1-\frac{3}{2}\lambda_1)\oplus
L(l\lambda_1-\frac{3}{2}\lambda_1),
\end{eqnarray*}
\end{enumerate}
\end{lemma}

Another way of realizing this decomposition is based on the metaplectic Howe duality, \cite{bss}. 
The metaplectic analogue of the classical theorem on the separation of variables allows to decompose 
the space $\mathrm{Pol}(\mR^2)\otimes {\mathcal S}$ of complex polynomials valued in the Segal-Shale-Weil
representation under the action of $\mp(2,\mR)$ into a direct sum of simple weight 
$\mp(2,\mR)$-modules
\begin{eqnarray}
\mathrm{Pol}(\mR^2)\otimes {\mathcal S}\simeq\bigoplus_{l=0}^\infty\bigoplus_{j=0}^\infty X_s^j{M}_l,
\end{eqnarray}
where we use the notation ${M}_l:={M}^+_l\oplus {M}^-_l$. This decomposition takes
the form of an infinite triangle 
\begin{eqnarray}\label{obrdekonposition}
\xymatrix@=11pt{\mathrm{P}_0 \otimes {\mathcal S} \ar@{=}[d] &  \mathrm{P}_1 \otimes {\mathcal S} \ar@{=}[d]& 
\mathrm{P}_2 \otimes {\mathcal S} \ar@{=}[d] & \mathrm{P}_3 \otimes {\mathcal S} \ar@{=}[d] & 
\mathrm{P}_4 \otimes {\mathcal S} \ar@{=}[d]& \mathrm{P}_5 \otimes {\mathcal S} \ar@{=}[d] \\
M_0 \ar[r] & X_s M_0 \ar @{} [d] |{\oplus} \ar[r] & X_s^2 M_0 \ar @{} [d] |{\oplus} \ar[r] & X_s^3 M_0 \ar @{} [d] |{\oplus}
 \ar[r] & X_s^4 M_0 \ar @{} [d] |{\oplus}\ar[r] & X_s^5 M_0 \ar @{} [d] |{\oplus} \\
& M_1 \ar[r] & X_s M_1 \ar @{} [d] |{\oplus}\ar[r] & X_s^2 M_1 \ar @{} [d] |{\oplus}
 \ar[r] & X_s^3 M_1 \ar @{} [d] |{\oplus}\ar[r] & X_s^4 M_1 \ar @{} [d] |{\oplus} \\
&& M_2 \ar[r] & X_s M_2 \ar @{} [d] |{\oplus}
 \ar[r] & X_s^2 M_2 \ar @{} [d] |{\oplus}\ar[r] & X_s^3 M_2 \ar @{} [d] |{\oplus} \\
&&& M_3 \ar[r] & X_s M_3 \ar @{} [d] |{\oplus}\ar[r] & X_s^2 M_3  \ar @{} [d] |{\oplus} \\
&&&& M_4 \ar[r] & X_s M_4 \ar @{} [d] |{\oplus}  \\
&&&&& M_5 
}
\end{eqnarray}

Let us now explain the notation used in the previous scheme. First of all, we used the shorthand 
notation $\mathrm{P}_l=\mathrm{Pol}_l(\mR^2),\, l\in\mN_0$, and all spaces and arrows on the picture have the following meaning.
The three operators ($i\in\mC$ is the complex unit)
\begin{eqnarray}
& & X_s={} y \partial_q + i x q,
\nonumber \\
& & D_s={} i q \partial_y - \partial_x\partial_q,
\nonumber \\
& & E=x\partial_x+y\partial_y, \label{xdeformulas}
\end{eqnarray}
where $D_s$ acts horizontally as $X_s$ but in the opposite direction, fulfill the $\sl(2,\mR)$-commutation relations:
\begin{eqnarray}
& & [E ,D_s]=-D_s,
\nonumber \\
\label{slRels}
& & [E ,X_s]=X_s,\\
& & [D_s,X_s]=-i(E +1).\nonumber 
\end{eqnarray}
Let $s(x,y,q)\in \mathrm{Pol}(\mR^2)\otimes {\mathcal S}$, $h\in \Mp(2,\mathbb{R})$ and 
$\pi (h)=g\in \Sp(2,\mathbb{R})$. We define
the action of $\Mp(2,\mR)$ to be
\begin{eqnarray}
& & \tilde{\varrho} (h)
s(x,y,q)=\varrho ( h) s(\pi(g^{-1}) \begin{pmatrix}
x  \\
y  
\end{pmatrix},q)=\varrho (h) s(dx-by,-cx+ay,q),
\nonumber \\
& & g=
\begin{pmatrix}
a & b \\
c & d 
\end{pmatrix}
\in \SL(2,\mR).
\end{eqnarray} 
where $\varrho$ acts on the Segal-Shale-Weil representation via (\ref{sshrepr}).
Passing to the infinitesimal action, we get the operators representing the basis elements of $\mp(2,\mR)$:
\begin{align*}
\frac{d}{dt}\Big|_{t=0}\tilde{\varrho} (\exp(tX))s(x,y,q)
={}&\frac{d}{dt}\Big|_{t=0}
\varrho \begin{pmatrix}
1 & t \\
0 & 1 
\end{pmatrix} 
s(x-yt,y,q)\\
={} &-\frac{i}{2} q^2 e^{-\frac{i}{2} t q^2}s(x-yt,y,q)\Big|_{t=0}\\
+{} & e^{-\frac{i}{2} t q^2} \frac{d}{dt}s(x-yt,y,q)\Big|_{t=0}\\
={} & \big(-\frac{i}{2} q^2 - y \frac{\partial}{\partial x}\big) s(x,y,q),
\end{align*}

\begin{align*}
\frac{d}{dt}\Big|_{t=0}\tilde{\varrho} (\exp(tH))
s(x,y,q)={}&\frac{d}{dt}\Big|_{t=0}
\varrho \begin{pmatrix}
e^t & t \\
0 & e^{-1} 
\end{pmatrix} 
s(x e^{-t},y e^{t},q)\\
={} &\frac{1}{2} e^{\frac{1}{2} t}
s(x e^{-t},y e^{t},qe^t) + e^{\frac{1}{2} t} \frac{d}{dt}s(x e^{-t},y e^{t},qe^t)\Big|_{t=0}\\
={} & \big(\frac{1}{2}- x \frac{\partial}{\partial x}+y \frac{\partial}{\partial y}+q \frac{\partial}{\partial q}\big) s(x,y,q),
\end{align*}

\begin{eqnarray}
& & \tilde{\varrho}(X)=- y \frac{\partial}{\partial x}-\frac{i}{2} q^2,\,
\tilde{\varrho}(Y)=- x \frac{\partial}{\partial y}-\frac{i}{2} \frac{\partial^2}{\partial q^2},\,
\nonumber \\
& & \tilde{\varrho}(H)=- x \frac{\partial}{\partial x}+y \frac{\partial}{\partial y}+q \frac{\partial}{\partial q}+\frac{1}{2},
\end{eqnarray}
and they satisfy commutation rules of $\mp(2,\mathbb{R})$:
\begin{align*}
[\tilde{\varrho}(X),\tilde{\varrho}(Y)]&=\tilde{\varrho}(H),\\
[\tilde{\varrho}(H),\tilde{\varrho}(X)]&=2\tilde{\varrho}(X),\\
[\tilde{\varrho}(H),\tilde{\varrho}(Y)]&=-2\tilde{\varrho}(Y).
\end{align*}
Notice that we have not derived the explicit formula for $\tilde{\varrho}(Y)$, because 
it easily follows from the Lie algebra structure. 
The action of the Casimir operator $\mathrm{Cas}\in \mathcal{U}(\mp(2,\mR))\otimes Cl_s(\mR^{2},\omega)$:
$$
\mathrm{Cas}=\tilde{\varrho}(H)^2+1+2\tilde{\varrho}(X)\tilde{\varrho}(Y)+2 \tilde{\varrho}(Y)\tilde{\varrho}(X),
$$
is given by the differential operator
\begin{eqnarray}
\mathrm{Cas} &=& x^2 \partial^2_x+ y^2 \partial^2_y +2x\partial_x +4y\partial_y + 2xy\partial_x\partial_y +\frac{1}{4}
\nonumber \\
& & -2xq\partial_x \partial_q +2 yq \partial_y\partial_q +2iy\partial_x\partial^2_q +2ix q^2\partial_y 
\nonumber \\
& = & E_x(E_x-1)+E_y(E_y-1)+2 E_x+4 E_y+2 E_x E_y+\frac{1}{4}
\nonumber \\
& & -2 E_x E_q +2 E_y E_q +2iy\partial_x\partial^2_q +2ix q^2\partial_y.
\end{eqnarray}
Here we introduced the notation $\partial_x:=\frac{\partial}{\partial x},\,\partial_x:=\frac{\partial}{\partial x}$
and $E_x=x\partial_x,\, E_y=y\partial_y,\, E_q=q\partial_q$ for the Euler homogeneity operators.
\begin{lemma}
The operators $X_s$ and $D_s$ commute with the operators $\tilde{\varrho}(X)$, $\tilde{\varrho}(Y)$ and $\tilde{\varrho}(H)$.
In other words, they are $\mp(2,\mR)$-invariant differential operators on complex polynomials valued in the Segal-Shale-Weil
representation. 
\end{lemma}
{\bf Proof:}

For example, we have
\begin{eqnarray}
[D_s,\tilde{\varrho}(H)]=i q \partial_y[\partial_y,y]+ i q\partial_q[q,\partial_q]+\partial_x\partial_q [\partial_x,x]-\partial_x\partial_q [\partial_q,q]=0,
\end{eqnarray}
and all remaining commutators are computed analogously.
\hfill $\square$

The action of $\mp(2,\mR)\times \sl(2,\mR)$ generates the multiplicity free decomposition
of the representation and the pair of Lie algebras in the product is called the metaplectic 
Howe dual pair. The operators $X_s,D_s$ act on the previous picture horizontally and 
isomorphically identify the neighboring $\mp(2,\mR)$-modules. The modules $M_l$, $l\in\mN$,
on the left-most diagonal are termed symplectic monogenics, and are characterized as $l$-homogeneous solutions of the symplectic Dirac 
operator $D_s$. Thus the decomposition is given as a tensor product of the symplectic monogenics
multiplied by polynomial algebra of invariants $\mC[X_s]$. The operator $X_s$ maps polynomial symplectic
spinors valued in the odd part 
of $\mathcal{S}$ into symplectic spinors valued in the even part of $\mathcal{S}$. This means that 
$M_m^-$ is valued in $\mathcal{S}_-$, $X_s M_m^-$ is valued in $\mathcal{S}_+$, etc.


\section{The symplectic twistor operator $T_s$}

We start with an abstract definition of the symplectic twistor operator $T_s$ 
and then we specialize to the standard symplectic space $(\mR^2,\omega)$.
\begin{definition}
Let $(M, \nabla ,\omega)$ be a symplectic spin manifold of dimension $2n$, $\nabla^s$ 
the associated symplectic spin covariant derivative
and $\omega\in C^\infty(M,\wedge^2T^\star M)$ a non-degenerate $2$-form such that $\nabla\omega=0$. We denote by 
$$
\{e_1,\ldots,e_{2n}\}\equiv \{e_1,\dots ,e_n,f_1,\dots ,f_n\}
$$ 
a local 
symplectic frame. The symplectic twistor operator $T_s$ on $M$ is the first order differential operator 
$T_s$ acting on smooth symplectic spinors ${\mathcal S}$:
\begin{eqnarray}
& & \nabla^s :\, C^\infty(M,{\mathcal S})\longrightarrow T^\star M\otimes 
C^\infty(M,{\mathcal S}),\nonumber \\
& & T_s:=P_{\mathrm{Ker}(c)}\circ\omega^{-1}\circ\nabla^s:\, C^\infty(M,{\mathcal S})
\longrightarrow C^\infty(M,{\mathcal T}),
\end{eqnarray}
where ${\mathcal T}$ is the space of symplectic twistors, 
$T^\star M\otimes{\mathcal S}\simeq{\mathcal S}\oplus{\mathcal T}$, given by algebraic projection 
$$
P_{\mathrm{Ker}(c_s)}:T^\star M\otimes C^\infty(M,{\mathcal S})\longrightarrow C^\infty(M,{\mathcal T})
$$ 
on the kernel of the symplectic Clifford multiplication $c_s$. In the local symplectic coframe 
$\{\epsilon^1\}^{2n}_{j=1}$ dual to the symplectic frame $\{e_j\}^{2n}_{j=1}$ with respect to $\omega$,
we have the local formula for $T_s$:
\begin{align}\label{localformula}
T_s=\Big(1+\frac{1}{n} \Big)\sum_{k=1}^{2n}\epsilon^k \otimes \nabla_{e_k}^s 
+ \frac{i}{n}\sum_{j,k,l=1}^{2n}\epsilon^l \otimes \omega^{kj} e_j \cdot e_l \cdot \nabla_{e_k}^s ,
\end{align} 
where $\cdot$ is the shorthand notation for the symplectic Clifford multiplication and $i\in\mC$ is the imaginary unit. 
We use the convention $\omega^{kj}=1$ for $j=k+n$ and $k=1,\dots ,n$, $\omega^{kj}=-1$ for $k=n+1,\dots ,2n$ 
and $j=k-n$, and $\omega^{kj}=0$ otherwise.
\end{definition}
The symplectic Dirac operator $D_s$, defined as the image of the symplectic Clifford multiplication $c_s$,
has the explicit form (\ref{xdeformulas}).
\begin{lemma}
The symplectic twistor operator $T_s$ is $\Mp(2n,\mR)$-invariant. 
\end{lemma}
{\bf Proof:}

The property of invariance is a direct consequence of the equivariance of symplectic 
covariant derivative and the invariance of algebraic projection $P_{\mathrm{Ker}(c_s)}$, and 
amounts to show that
\begin{equation}\label{invar}
T_s( \tilde{\varrho}(g)s)=\pi(g)\otimes\tilde{\varrho}(g)( T_s s)
\end{equation}
for any $g\in \Mp(2n,\mR)$ and $s\in C^\infty(M,{\mathcal S})$.
Using the local formula (\ref{localformula}) for $T_s$ in a local chart 
$(x_1,\ldots,x_{2n})$, both sides of (\ref{invar}) are equal  
\begin{align*}
{}&\Big(1+\frac{1}{n} \Big)\sum_{k=1}^{2n}\epsilon^k \otimes \varrho(g) 
\frac{\partial}{\partial x_k} \big[s\big(\pi(g)^{-1}x\big)\big]\\
&+ \frac{i}{n}\sum_{j,k,l=1}^{2n}\epsilon^l \otimes \omega^{kj} e_j \cdot e_l \cdot 
\Big[\varrho(g) \frac{\partial}{\partial x_k} \big[s\big(\pi(g)^{-1}x\big)\big]\Big]
\end{align*}
and the proof follows.

\hfill
$\square$

In the case $M=(\mathbb{R}^{2n},\omega)$, the symplectic twistor operator is  
\begin{equation}\label{TvR}
T_s=\Big(1+\frac{1}{n} \Big)\sum_{k=1}^{2n}\epsilon^k \otimes \frac{\partial}{\partial x_k}  
+ \frac{i}{n}\sum_{j,k,l=1}^{2n}\epsilon^l \otimes \omega^{kj} e_j \cdot e_l \cdot \frac{\partial}{\partial x_k}.
\end{equation}
\begin{lemma}\label{tvar}
In the case of the standard symplectic space $(\mR^2,\omega)$  with coordinates $x,y$ and $\omega=dx\wedge dy$, a symplectic frame $\{e, f\}$
and its dual coframe $\{\epsilon^1,\epsilon^2\}$,
the symplectic twistor operator 
$T_s: C^\infty(\mathbb{R}^2,\mathcal{S})\rightarrow$ $ C^\infty(\mathbb{R}^2, {\mathcal T})$
acts on a smooth symplectic spinor $s(x,y,q)\in {C}^\infty(\mathbb{R}^2,{\mathcal S})$ by
\begin{equation}\label{tscomp}
T_s (s)=\epsilon^1 \otimes\Big(\frac{\partial s}{\partial x}- q\frac{\partial^2 s}{\partial q \partial x}+iq^2\frac{\partial s}{\partial y} \Big)
+\epsilon^2 \otimes\Big(2 \frac{\partial s}{\partial y} +i\frac{\partial^3 s}{\partial q^2 \partial x}+q \frac{\partial^2 s}{\partial q\partial y} \Big).
\end{equation}
\end{lemma}
The last display follows from (\ref{TvR}) by direct substitution for the symplectic Clifford multiplication. 
The next Lemma simplifies the condition on a symplectic spinor to be in the kernel of $T_s$.
\begin{lemma}\label{halftwistor}
A smooth symplectic spinor $s(x,y,q)\in {C}^\infty(\mathbb{R}^2,{\mathcal S})$
is in the kernel of $T_s$ if and only if it fulfils the partial differential equation 
\begin{eqnarray}\label{rce}
\Big(\frac{\partial}{\partial x}-q \frac{\partial^2}{\partial q \partial x}+ i q^2 \dfrac{\partial}{\partial y}\Big)s=0.
\end{eqnarray} 
\end{lemma}
{\bf Proof:} 

The claim is a consequence of Lemma \ref{tvar}, because the covectors $\epsilon^1, \epsilon^2$ are linearly independent
and the differential operators in (\ref{tscomp}) (the two components of $T_s$ by 
$\epsilon^1$ and $\epsilon^2$) have the same solution space
(i.e., $s$ solving one of them implies that $s$ solves the second one.) This implies the equivalence statement in the
Lemma.

\hfill $\square$

Notice that $\tilde{\varrho}(X)$, $\tilde{\varrho}(Y)$ and $\tilde{\varrho}(H)$ preserve the solution space 
of the twistor equation (\ref{rce}), i.e.,
if the symplectic spinor $s$ solves (\ref{rce}) then $\tilde{\varrho}(X)s,\tilde{\varrho}(Y)s$ and $\tilde{\varrho}(H)s$ solve (\ref{rce}) . 
This is a consequence of $\Mp(2,\mR)$-invariance of the symplectic twistor operator $T_s$ on $\mR^2$ (in fact, the same observation is true in any dimension.)
By abuse of notation, we use $T_s$ in Section $3$ to denote the operator (\ref{rce}) and call it 
the symplectic twistor operator - this terminology is justified by the reduction in Lemma \ref{halftwistor}.
In the article, we work with polynomial (in $x,y$ or $z,\bar z$) smooth symplectic spinors $\mathrm{Pol}(\mathbb{R}^2,{\mathcal S})$.


\section{The polynomial solution space of the symplectic twistor operator $T_s$ on $\mR^2$}

Let us consider the complex vector space of symplectic spinor valued polynomials 
$\mathrm{Pol}(\mathbb{R}^2,{\mathcal S})$, 
${\mathcal S}\simeq {\mathcal S}_-\oplus {\mathcal S}_+$, together with its decomposition
on irreducible subspaces with respect to the natural action of $\mp(2,\mR)$. It follows 
from the $\mp(2,\mR)$-invariance of the symplectic twistor operator that it is sufficient to 
characterize its behavior on any non-zero vector in an irreducible $\mp(2,\mR)$-submodule,
and that its action preserves the subspace of homogeneous symplectic spinors.
This is what we are going to accomplish in the present section. Note that the meaning of the 
natural number $n\in\mN$ used in previous sections to denote the dimension of the 
underlying symplectic space is different from its use in the present section.    

The main technical difficulty consists of finding suitable representative smooth vectors 
in each irreducible $\mp(2,\mR)$-subspace. We shall find a general characterizing 
condition for a polynomial (in the real variables $x,y$) valued in the Schwartz space $S(\mR)$ (in the 
variable $q$) as a formal power series,
and the representative vectors are always conveniently chosen as polynomials in $q$ weighted by the exponential $e^{-\frac{q^2}{2}}$. 
In other words, the constructed vectors are $\tilde K=\widetilde{\mathrm{U}}(1)$-finite vectors in 
$S(\mR)$.
These representative vectors are then evaluated on the symplectic twistor operator $T_s$.

First of all, the constant symplectic spinors belong to the solution space of $T_s$. 
We have
\begin{lemma}
\begin{align}
T_s (X_s e^{-\frac{q^2}{2}})=\,&T_s(i e^{-\frac{q^2}{2}} q (x + i y))=0,\\
T_s (X_s qe^{-\frac{q^2}{2}})=\,&T_s(e^{-\frac{q^2}{2}} (i q^2 (x + i y) + y))=0.
\end{align}
\end{lemma}
The next Lemma is preparatory for further considerations.
\begin{lemma}\label{Xs^n}
We have for any $n\in\mN_0$, $(X_s)^n\in \End(\mathrm{Pol}(\mathbb{R}^2,{\mathcal S}))$, the following identity
\begin{equation}\label{XsMe}
(X_s)^n=\sum^{\left \lfloor\frac{n}{2}\right \rfloor}_{j=0} \sum^{n-2j}_{k=0}A_{jk}^n y^{n-j-k}(ix)^{j+k}q^k\partial_q^{n-2j-k}.
\end{equation}
Here $\left \lfloor\frac{n}{2}\right \rfloor$ is the floor function applied to $\frac{n}{2}$, and the coefficients 
$A_{jk}^n\in\mC$ fulfill the $4$-term recurrent relation
\begin{eqnarray}\label{recur}
A_{jk}^n=A_{jk}^{(n-1)}+A_{j(k-1)}^{(n-1)}+(k+1)A_{(j-1)(k+1)}^{(n-1)}.
\end{eqnarray}
We use the normalization $A_{00}^0=1$, and $A_{jk}^n\neq 0$ only for $n\in\mN_0$, 
$j=0,\ldots,\left \lfloor\frac{n}{2}\right \rfloor$, and $k=0,\ldots,n-2j$.
\end{lemma}
{\bf Proof:}

The proof is by induction on $n\in\mN_0$. The claim is trivial for $n=0$, 
and for $n=1$ we have 
$$
(X_s)^1=A_{00}^1 y\partial_q+A_{01}^1 ix q,
$$ 
where $A_{00}^1=A_{00}^0=1$ and $A_{01}^1=A_{00}^0=1$.

We assume that the formula holds for $n-1$ and aim to prove it for $n$:
{\footnotesize
\begin{align*}
&{}(ixq + y\partial_q)\Big(\sum^{\left \lfloor\frac{n-1}{2}\right \rfloor}_{j=0} \sum^{n-1-2j}_{k=0}A_{jk}^{(n-1)}y^{n-1-j-k}(ix)^{j+k}q^k\partial_q^{n-1-2j-k} \Big)\\
=&{} \sum^{\left \lfloor\frac{n-1}{2}\right \rfloor}_{j=0} \sum^{n-1-2j}_{k=0}A_{jk}^{(n-1)} \Big(y^{n-1-j-k}(ix)^{j+k+1}q^{k+1}\partial_q^{n-1-2j-k} \\
&{}+y^{n-j-k}(ix)^{j+k}q^k\partial_q^{n-2j-k}  + k y^{n-j-k}(ix)^{j+k}q^{k-1}\partial_q^{n-1-2j-k}\Big)\\
=&{}\sum^{\left \lfloor\frac{n-1}{2}\right \rfloor}_{j=0} \sum^{n-2j}_{k=0}A_{j(k-1)}^{(n-1)}y^{n-j-k}(ix)^{j+k}q^k\partial_q^{n-2j-k}\\
&{}+\sum^{\left \lfloor\frac{n-1}{2}\right \rfloor}_{j=0} \sum^{n-2j}_{k=0}A_{jk}^{(n-1)}y^{n-j-k}(ix)^{j+k}q^k\partial_q^{n-2j-k}\\
&{}+\sum^{\left \lfloor\frac{n-1}{2}\right \rfloor+1}_{j=0} \sum^{n-2j}_{k=0}(k+1)A_{(j-1)(k+1)}^{(n-1)}y^{n-j-k}(ix)^{j+k}q^k\partial_q^{n-2j-k}
\end{align*}}
where we shifted the indexes in the first sum by $k\rightarrow k-1$, in the third sum by $k\rightarrow k+1$ and $j\rightarrow j-1$. Altogether we get.
{\footnotesize
\begin{align*}
& \!\sum^{\left \lfloor\frac{n-1}{2}\right \rfloor}_{j=0} \sum^{n-2j}_{k=0}\!\!\!\Big(\!
A_{jk}^{(n-1)}\!+\! A_{j(k-1)}^{(n-1)}\!+\!(k+1)A_{(j-1)(k+1)}^{(n-1)} \Big)y^{n-j-k}(ix)^{j+k}q^k\partial_q^{n-2j-k}\\
& \!\!+\!\!\!\!\sum^{n-2\left \lfloor\frac{n-1}{2}\right \rfloor-2}_{k=0}\!\!\!\!\!\!(k+1)A_{\left \lfloor\frac{n-1}{2}\right \rfloor(k+1)}^{(n-1)}\!y^{n-\left \lfloor\frac{n-1}{2}\right \rfloor-1-k}(ix)^{\left \lfloor\frac{n-1}{2}\right \rfloor+1+k}q^k\partial_q^{n-2\left \lfloor\frac{n-1}{2}\right \rfloor-2-k}
\end{align*}}
Now we apply the induction argument to the first term. The second term is non zero only for even $n$, when the previous expression equals to
$$
\sum^{\left \lfloor\frac{n}{2}\right \rfloor}_{j=0} \sum^{n-2j}_{k=0}A_{jk}^{n}y^{n-j-k}(ix)^{j+k}q^k\partial_q^{n-2j-k},
$$
which completes the required statement.

\hfill $\square$

\begin{remark}
Notice that for $j=0$, the solution of recurrent relation in (\ref{recur}) corresponds to the binomial coefficients. 
It follows from $A_{(-1)(k+1)}^{(n-1)}=0$,
$$
A_{0k}^n=A_{0k}^{(n-1)}+A_{0(k-1)}^{(n-1)},
$$
and therefore, $A_{0k}^n=\binom{n}{k}$.
\end{remark}
\begin{lemma}
We have 
$A_{1 (n-2)}^n=\frac{n(n-1)}{2}=\binom{n}{n-2}$. 
\end{lemma}
{\bf Proof:}

We use the relation 
$A_{1 (n-2)}^n=A_{1 (n-2)}^{(n-1)}+A_{1 (n-3)}^{(n-1)}+(n-1)A_{0 (n-1)}^{(n-1)}$,
where $A_{1 (n-2)}^{n-1}=0$ (because it is out of the range for the index $k$ in the equation (\ref{recur}).)
The proof goes by induction in $n$: we start with $A_{10}^2=A_{01}^1=1$, and claim $A_{1 (n-2)}^n=\frac{n(n-1)}{2}$. 
The induction step gives $A_{1 (n-1)}^{(n+1)}=A_{1 (n-2)}^{n}+n A_{0 n}^n=\frac{n^2-n}{2}+n=\frac{n^2+n}{2}.$

\hfill 
$\square$

A direct consequence of the Baker-Campbell-Hausdorff formula or 
its dual, Zassenhaus formula, for
three operators $A, B, C$ fulfilling the commutation relations $[A,B]=C$ and 
$[A,C]=[B,C]=0$:
\begin{equation}
(A+B)^n=\sum_{\substack{
   k\leq n\\
   k \equiv n \mbox{ (mod }2 \mbox{)}}}\left(\sum_{r=0}^k \binom{k}{r}A^r B^{k-r}  \right) \left(-\frac{C}{2} \right)^{\frac{n-k}{2}} \frac{n!}{k!\left( \frac{n-k}{2} \right)!},
\end{equation}
is the solution of the recursion relation \eqref{recur}.
\begin{theorem}
For $X_s= i x q+y \partial_q$ with $[ i x q, y \partial_q]= -ixy$, we have
\begin{equation}\label{XsBCH}
(i x q+y \partial_q)^n=\sum_{\substack{
   k\leq n\\
   k \equiv n\mbox{ (mod }2 \mbox{)}}} \left(\sum_{r=0}^k \frac{n!}{r!(k-r)!\left( \frac{n-k}{2} \right)!2^{\frac{n-k}{2}}} (i x)^{r+\frac{n-k}{2}} y^{k-r+\frac{n-k}{2}}  q^r  \partial_q^{k-r} \right)
\end{equation}
for any $n\in\mN_0$, and the comparison with 
Lemma \ref{Xs^n} yields the solution of the recursion 
relation \eqref{recur}:
\begin{equation}
A^n_{jk}=\frac{n!}{k!(n-2j-k)!j!2^j},
\end{equation}
where the index $k$ in \eqref{XsBCH} corresponds to $n-2j$ in \eqref{XsMe}, 
and the index $r$ in \eqref{XsBCH} corresponds to $k$ in \eqref{XsMe}.
\end{theorem}

Let us remark that the composition $T_s\circ (X_s)^n$ for $n=2,3$, acting on $e^{-\frac{q^2}{2}}$ and $qe^{-\frac{q^2}{2}}$, 
is non-vanishing. This means that some irreducible $\mp(2,\mR)$-components in the decomposition (\ref{obrdekonposition}) are not in the kernel of $T_s$:

\begin{align}
T_s(X_s^2 e^{-\frac{q^2}{2}})=e^{-\frac{q^2}{2}}&{} (q^2 x + i y + i q^2 y)\not=0, \nonumber\\ \nonumber
T_s(X_s^2 q e^{-\frac{q^2}{2}})=e^{-\frac{q^2}{2}}&{} (q^3 x + i q^3 y)\not=0,\\ \nonumber
T_s(X_s^3 e^{-\frac{q^2}{2}})= e^{-\frac{q^2}{2}}&{}(3 i q^3 x^2 - 6 q^3 x y - 3 i q^3 y^2)\not=0,\\ \nonumber
T_s(X_s^3 q e^{-\frac{q^2}{2}})=e^{-\frac{q^2}{2}}&{} (3i q^4 x^2 +6 q^2 x y -6 q^4 xy + 3i y^2 + 6iq^2 y^2\\
&{} -3i q^4 y^2)\not=0, 
\end{align}

\begin{lemma}
Let $n\in\mN_0$. Then
{\footnotesize
\begin{align}\label{TXs^n}
T_s&{}\circ (X_s)^n=\sum^{\left \lfloor\frac{n}{2}\right \rfloor}_{j=0} \sum^{n-2j}_{k=0}A_{jk}^n\Big(
 i(j+k) y^{n-j-k}(ix)^{j+k-1}q^k\partial_q^{n-2j-k}\nonumber\\
&{}+ y^{n-j-k}(ix)^{j+k}q^k\partial_x\partial_q^{n-2j-k} 
 -i(j+k) y^{n-j-k}(ix)^{j+k-1}q^{k+1}\partial_q^{n-2j-k+1}\nonumber\\
&{} - y^{n-j-k}(ix)^{j+k}q^{k+1}\partial_x\partial_q^{n-2j-k+1}
 -ik(j+k) y^{n-j-k}(ix)^{j+k-1}q^k\partial_q^{n-2j-k}\nonumber\\
&{} -k y^{n-j-k}(ix)^{j+k}q^k\partial_x\partial_q^{n-2j-k}
  \!+i(n-j-k) y^{n-j-k-1}(ix)^{j+k}q^{k+2}\partial_q^{n-2j-k}\nonumber\\
&{}+ i y^{n-j-k}(ix)^{j+k}q^{k+2}\partial_y\partial_q^{n-2j-k}
\Big)
\end{align}}
In particular, 
$T_s((X_s)^n e^{-\frac{q^2}{2}}) \neq 0$ and $T_s((X_s)^n qe^{-\frac{q^2}{2}}) \neq 0$ for all $n>1$.
\end{lemma}
{\bf Proof:}

The proof is based on the identity in Lemma \ref{Xs^n}.
The non-triviality of the composition is detected by the coefficient in the monomial 
$x^{n-1} q^{n} e^{-\frac{q^2}{2}}$ in $T_s((X_s)^n e^{-\frac{q^2}{2}})$. 
It follows from the identity (\ref{TXs^n}) that this coefficient is 
\begin{align}
&{}i^n(A_{0n}^n n- A_{0n}^n n^2 +A_{1 (n-2)}^n)x^{n-1} q^{n} e^{-\frac{q^2}{2}}\nonumber \\
=&{}i^n\left(\binom{n}{n}(n-  n^2)+\binom{n}{n-2} \right)x^{n-1} q^{n} e^{-\frac{q^2}{2}}\nonumber\\
=&{}-i^n\frac{n(n-1)}{2}  x^{n-1} q^{n} e^{-\frac{q^2}{2}},
\end{align}
which is non-zero for all $n>1$.

As for the action on the vector $qe^{-\frac{q^2}{2}}$, the situation is analogous. 
The coefficient of the monomial $x^{n-1} q^{n+1} e^{-\frac{q^2}{2}}$ in 
$T_s((X_s)^n qe^{-\frac{q^2}{2}})$ is $-i^n\frac{n(n-1)}{2}$, which is again non-zero for all $n>1$.
The proof is complete.

\hfill $\square$

In the next part we focus for a while on symplectic spinors given by iterative action of $X_s$ on 
${\mathcal S}_+$, and complete the task of finding all subspaces of polynomial solutions of 
$T_s$ (expressed in the real variables $x,y$).
\begin{lemma}\label{allevenmonogenics}
The vectors $e^{-\frac{q^2}{2}}(x + i y)^m\in \mathrm{Pol}_m(\mathbb{R}^2,{\mathcal S}_+)$, $m\in\mN_0$,   
are in the kernel of $D_s$, but not in the kernel of the symplectic twistor operator $T_s$.
\end{lemma}
{\bf Proof:}

 We get by direct computation,
\begin{align}
D_s\big(e^{-\frac{q^2}{2}}(x + i y)^m\big)=&{} iq\partial_y e^{-\frac{q^2}{2}}(x + i y)^m- \partial_x\partial_qe^{-\frac{q^2}{2}}(x + i y)^m
\nonumber \\ \nonumber
 =&{} e^{-\frac{q^2}{2}}(-m q (x+iy)^{m-1}+mq (x+iy)^{m-1})=0,
\nonumber \\ \nonumber
T_s\big(e^{-\frac{q^2}{2}}(x + i y)^m\big)=&{} \partial_x e^{-\frac{q^2}{2}} (x + i y)^m-q \partial_x\partial_q e^{-\frac{q^2}{2}} (x + i y)^m
\nonumber \\ \nonumber 
&{} +i q^2\partial_y e^{-\frac{q^2}{2}} (x + i y)^m=e^{-\frac{q^2}{2}} m(x + i y)^{m-1}\neq 0
\end{align}
for any natural number $m>0$.
\hfill $\square$

\begin{lemma}
Let $m\in\mN_0$.
Then the vectors $X_se^{-\frac{q^2}{2}}(x + i y)^m$ in $\mathrm{Pol}_{m+1}(\mathbb{R}^2,{\mathcal S}_+)$ 
are in the kernel of the symplectic twistor operator $T_s$.
\end{lemma}
{\bf Proof:}

We have
\begin{eqnarray}
& & T_s(X_s e^{-\frac{q^2}{2}}(x + i y)^m) = T_s(i q e^{-\frac{q^2}{2}}(x+i y)^{m+1})
\nonumber \\ \nonumber
& & =i (m+1)e^{-\frac{q^2}{2}}(q-q+q^2-q^2)(x+iy)^m=0.
\end{eqnarray}
\hfill 
$\square$

\begin{remark}
The non-trivial elements in $\mathrm{Ker}(T_s)$ are 
\begin{equation}
q e^{-\frac{q^2}{2}}(x+i y)^{k},\, k\in \mathbb{N}_0.
\end{equation}
\end{remark}

The next Lemma completes the information on 
the behaviour of $T_s$ for remaining $\mp(2,\mR)$-modules coming from the action of
$X_s$ on ${\mathcal S}_+$.
\begin{lemma}
For all natural numbers $n>1$ and all $m\in\mN_0$, we have 
\begin{eqnarray}
T_s((X_s)^n e^{-\frac{q^2}{2}}(x+iy)^m) \neq 0.
\end{eqnarray} 
\end{lemma}
{\bf Proof:}

We focus on the coefficient by the monomial $x^{n-1+m} q^{n} e^{-\frac{q^2}{2}}$ in $T_s((X_s)^n e^{-\frac{q^2}{2}})$. 
It follows from (\ref{TXs^n}) that the contribution to this coefficient is 
\begin{align}
i^n(A_{0n}^n n-&{} A_{0n}^n n^2 + A_{1 (n-2)}^n+A_{0 n}^n m-A_{0 n}^n m n)x^{n-1+m} q^{n} e^{-\frac{q^2}{2}}\nonumber\\
=&{}i^n\left(\binom{n}{n}(n-n^2+m-m n)+\binom{n}{n-2}\right)x^{n-1+m} q^{n} e^{-\frac{q^2}{2}}\nonumber\\
=&{}-i^n\frac{(n+2m)(n-1)}{2}  x^{n-1+m} q^{n} e^{-\frac{q^2}{2}},
\end{align}
which is non-zero for all natural numbers $n>1$ and all $m\in\mN_0$.

\hfill $\square$

Let us summarize the previous lemmas in the final Theorem.
\begin{theorem}\label{stevenelement}
The solution space of the symplectic twistor operator $T_s$ acting on $\mathrm{Pol}(\mathbb{R}^2,{\mathcal S}_+)$ consists
of the set of $\mp(2,\mR)$-modules in the boxes, realized in the decomposition of
$\mathrm{Pol}(\mathbb{R}^2,{\mathcal S}_+)$ on $\mp(2,\mR)$ irreducible subspaces:
\begin{eqnarray}
\xymatrix@=11pt{
*+[F]{M^+_0}\ar @{}[d] |{ e^{-\frac{q^2}{2}}} \ar[r] & *+[F]{X_s M^+_0} \ar @{} [d] |{\oplus} \ar[r] & X_s^2 M^+_0 \ar @{} [d] |{\oplus} \ar[r] & X_s^3 M^+_0 \ar @{} [d] |{\oplus}
 \ar[r] & X_s^4 M^+_0 \ar @{} [d] |{\oplus}\ar[r] & X_s^5 M^+_0 \ar @{} [d] |{\oplus} &\ldots\\
& M^+_1 \ar @{}[d] |{e^{-\frac{q^2}{2}}(x + i y)} \ar[r] & *+[F]{X_s M^+_1} \ar @{} [d] |{\oplus}\ar[r] & X_s^2 M^+_1 \ar @{} [d] |{\oplus}
 \ar[r] & X_s^3 M^+_1 \ar @{} [d] |{\oplus}\ar[r] & X_s^4 M^+_1 \ar @{} [d] |{\oplus} &\ldots\\
&& M^+_2\ar @{}[d] |{ e^{-\frac{q^2}{2}}(x + i y)^2} \ar[r] & *+[F]{X_s M^+_2} \ar @{} [d] |{\oplus}
 \ar[r] & X_s^2 M^+_2 \ar @{} [d] |{\oplus}\ar[r] & X_s^3 M^+_2 \ar @{} [d] |{\oplus} &\ldots\\
&&& M^+_3\ar @{}[d] |{ e^{-\frac{q^2}{2}}(x + i y)^3} \ar[r] & *+[F]{X_s M^+_3} \ar @{} [d] |{\oplus}\ar[r] & X_s^2 M^+_3  \ar @{} [d] |{\oplus}&\ldots\\
&&&& M^+_4\ar @{}[d] |{ e^{-\frac{q^2}{2}}(x + i y)^4} \ar[r] & *+[F]{X_s M^+_4} \ar @{} [d] |{\oplus}  &\ldots\\
&&&&& M^+_5 &\ldots
}
\end{eqnarray}
Notice that non-zero representative vectors in the solution space of $D_s$ are pictured under the spaces
of symplectic monogenics.
\end{theorem}

This completes the picture in the case of ${\mathcal S}_+$.
As we shall see, the representative solutions of $D_s$ for arbitrary homogeneity are far more complicated 
for ${\mathcal S}_-$ than for ${\mathcal S}_+$, which were chosen to be the powers of $z=x+iy$. A rather 
convenient way to simplify the presentation is to pass from the real coordinates $x,y$ to the complex coordinates
$z,\overline{z}$ for the standard complex structure on $\mR^2$, where 
$\partial_x=(\partial_z+\partial_{\bar{z}})$ and 
$\partial_y=i(\partial_z-\partial_{\bar{z}})$. 
\begin{lemma}
The operators $X_s,D_s$ and $T_s$
are in the complex coordinates $z, \bar{z}$ given by 
\begin{eqnarray}
& & X_s=\frac{i}{2}\big( (q-\partial_q)z+(q+\partial_q)\bar{z}\big),\label{XsHolom}
\nonumber \\
& & D_s=-\big((q+\partial_q)\partial_z+(-q +\partial_q)\partial_{\bar{z}} \big),
\label{DsHolom} \nonumber
\\
\label{THolom}
& & T_s=(1-q\partial_q-q^2)\partial_z+(1-q\partial_q+q^2)\partial_{\bar{z}}.
\end{eqnarray}
\end{lemma}
In the rest of the article we suppress the overall constants in 
$X_s,D_s,T_s$. The reason is that both the metaplectic Howe duality and the solution space of 
$D_s$, $T_s$ are independent of the normalization of $X_s$, $D_s$, $T_s$. In other words, the representative 
solutions differ by a 
non-zero multiple, a property which has no effect on our results. An element of the solution space of $D_s$ is called the symplectic monogenic.  

We start with the characterization of elements in the solution space of $D_s$, both for ${\mathcal S}_+$ and ${\mathcal S}_-$.

\begin{theorem}\label{DsRecurTheorem}
\begin{enumerate}
\item
The symplectic spinor of the homogeneity $m\in\mN_0$ in the variables $z,\bar{z}$,
\begin{equation}\label{DsRecurVyraz}
s=e^{-\frac{q^2}{2}}q\big(A^m(q)z^m+ A^{m-1}(q)z^{m-1}\bar{z}+ \ldots + A^1(q)z\bar{z}^{m-1}+A^0(q)\bar{z}^m \big),
\end{equation}
with coefficients in the formal power series in $q$,
$$
A^r(q)=a^r_0+a^r_2 q^2+a^r_4 q^4 + \ldots ,\, a^r_{k}\in\mC,\, r=0,\ldots ,m,\, k\in2\mN_0
$$
is in the kernel of $D_s$ provided the coefficients $a^r_k$ satisfy the system of recurrence relations
\begin{eqnarray}
& & 0=m(k+1)a^m_k+(k+1)a^{m-1}_k-2a^{m-1}_{k-2},\nonumber \\
& & 0=(m-1)(k+1)a^{m-1}_k+2(k+1)a^{m-2}_k-4a^{m-2}_{k-2},\nonumber \\
& & \,\,\,\dots \nonumber \\
& & 0=2(k+1)a^{2}_k+(m-1)(k+1)a^{1}_k-2(m-1)a^{1}_{k-2},\nonumber \\
& & 0=(k+1)a^1_k+m(k+1)a^{0}_k-2ma^{0}_{k-2},
\end{eqnarray}
equivalent to
\begin{equation}\label{DsRecurZkraceno}
(m-p)(k+1)a_k^{m-p}+(p+1)(k+1)a_k^{m-1-p}-2(p+1)a_{k-2}^{m-1-p}=0,
\end{equation}
for all $p=0,1,\ldots,m-1$.
\item
The symplectic spinor of the homogeneity $m\in\mN_0$ in the variables $z,\bar{z}$,
\begin{equation}
s=e^{-\frac{q^2}{2}}\big(A^m(q)z^m+ A^{m-1}(q)z^{m-1}\bar{z}+ \ldots + A^1(q)z\bar{z}^{m-1}+A^0(q)\bar{z}^m \big),
\end{equation}
with coefficients in the formal power series in $q$,
$$
A^r(q)=a^r_0+a^r_2 q^2+a^r_4 q^4 + \ldots ,\, a^r_{k}\in\mC,\, r=0,\ldots ,m,\, k\in2\mN_0
$$
is in the kernel of $D_s$ provided the coefficients $a^r_k$ satisfy the system of recurrence relations
\begin{eqnarray}
& & 0=m k a^m_k+k a^{m-1}_k-2a^{m-1}_{k-2}, \nonumber\\
& & 0=(m-1)k a^{m-1}_k+2 k a^{m-2}_k-4a^{m-2}_{k-2}, \nonumber\\
& & \,\,\, \dots \nonumber \\
& & 0=2 k a^{2}_k+(m-1) k a^{1}_k-2(m-1)a^{1}_{k-2},\\ \nonumber
& & 0= k a^1_k+m k a^{0}_k-2ma^{0}_{k-2},
\end{eqnarray}
equivalent to
\begin{equation}\label{DsRecurZkracSude}
(m-p)k a_k^{m-p}+(p+1) k a_k^{m-1-p}-2(p+1)a_{k-2}^{m-1-p}=0,
\end{equation}
for all $p=0,1,\ldots,m-1$.
\end{enumerate}
\end{theorem}
{\bf Proof:}

Because
$$
(q+\partial_q)e^{-\frac{q^2}{2}}qA^r(q)=e^{-\frac{q^2}{2}}[q^2+1-q^2+q\partial_q]A^r(q),
$$
$$
(-q+\partial_q)e^{-\frac{q^2}{2}}qA^r(q)=e^{-\frac{q^2}{2}}[-q^2+1-q^2+q\partial_q]A^r(q),
$$
the action of $D_s$ on the vector $e^{-\frac{q^2}{2}}qA^r(q)$ is
\begin{align}\label{DsRecurentDk}
D_s\Big(e^{-\frac{q^2}{2}}q\big(&{}A^m(q)z^m+ A^{m-1}(q)z^{m-1}\bar{z}+ \ldots + A^1(q)z\bar{z}^{m-1}+A^0(q)\bar{z}^m \big)
\Big)\nonumber\\
=e^{-\frac{q^2}{2}}\Big(&{}
z^{m-1}\big(m[1+q\partial_q]A^m(q)+[1+q\partial_q-2q^2]A^{m-1}(q)\big)\nonumber\\
&{}z^{m-2}\bar{z}\big((m-1)[1+q\partial_q]A^{m-1}(q)+2[1+q\partial_q-2q^2]A^{m-2}(q)\big)\nonumber\\
&{}\vdots\nonumber\\
&{}z\bar{z}^{m-1}\big(2[1+q\partial_q]A^2(q)+(m-1)[1+q\partial_q-2q^2]A^1(q)\big)\nonumber\\
&{}\bar{z}^m\big([1+q\partial_q]A^{1}(q)+m[1+q\partial_q-2q^2]A^0(q)\big)
\Big).
\end{align}
The action of $[1+q\partial_q]$ on $A^r(q)$ yields 
$\sum_{k\in 2\mN}(k+1)a^r_k q^k$, and the action of 
$[1+q\partial_q-2q^2]$ on $A^r(q)$ gives $\sum_{k\in 2\mN}((k+1)a^r_k-2a^r_{k-2}) q^k$, for all $r=0,\ldots ,m$.

As for the second part, we have 
\begin{align*}
(q+\partial_q)e^{-\frac{q^2}{2}}A^r(q)=\,& e^{-\frac{q^2}{2}}[\partial_q]A^r(q),\\
(-q+\partial_q)e^{-\frac{q^2}{2}}A^r(q)=\,& e^{-\frac{q^2}{2}}[-2q+\partial_q]A^r(q),
\end{align*}
and the rest of the proof is analogous to the first part.
The proof is complete.

\hfill $\square$

\begin{remark}
\label{remarkexplicitmonogenic}
We observe that the choice of the constant $A^0(q)=a_0^0\neq 0$, i.e. $a^0_k\neq 0$ only for $k=0$, 
leads to the solution (polynomial in $q$) of the recurrence relation for all coefficients in the symplectic
spinor (\ref{DsRecurVyraz}). 
\begin{align*}
A^0(q)=&{}a_0^0,\\
A^1(q)=&{}\left(-1 +\frac{2}{3}q^2\right)\binom{m}{1} a_0^0,\\
&{}\ldots\\
A^r(q)=&{}\left((-1)^r+\ldots +\frac{2^r }{(2r+1)!!}q^{2r} \right)\binom{m}{r}a_0^0,\\
&{}\ldots\\
A^m(q)=&{}\left((-1)^m + \ldots +\frac{2^m }{(2m+1)!!}q^{2m} \right)\binom{m}{m}a_0^0, 
\end{align*}
where $(2m+1)!!=(2m+1)\cdot(2m-1)\cdots 3\cdot 1$.
In this way, we get simple representative vectors in the kernel of $D_s$,
valued in ${\mathcal S}_-$ for each homogeneity $m$. We have for $m=1,2,3$
\begin{align}
e^{-\frac{q^2}{2}}q\bigg(&{}\Big(-1 + \frac{2}{3} q^2\Big)z  + \bar{z}\bigg) a_0^0, \nonumber\\ \nonumber
e^{-\frac{q^2}{2}}\bigg(q&{} \Big(1 - \frac{4}{3} q^2 + \frac{4}{15} q^4\Big)z^2 +  \Big(-2 + \frac{4}{3} q^2\Big)z \bar{z} + \bar{z}^2\bigg)a_0^0,\label{DsPrikHom2}\\
e^{-\frac{q^2}{2}}\bigg(q&{} \Big(-1 + 2 q^2 - \frac{12}{15} q^4 + \frac{8}{105} q^6\Big)z^3 + 
   \Big(3 - 4 q^2 + \frac{4}{5} q^4 \Big)z^2 \bar{z} \nonumber\\
    +&{} (-3 + 2 q^2)z \bar{z}^2 + \bar{z}^3\bigg)a_0^0.
\end{align}
The same formulas expressed in the real variables $x,y$:
\begin{align}
\frac{2}{3}e^{-\frac{q^2}{2}}&{}\Big(q^3 (x + i y) - 3 i q y\Big)a_0^0,\nonumber \\
\frac{4}{15}e^{-\frac{q^2}{2}}&{}\Big( q^5 (x + i y)^2  + 10 q^3 y (-i x + y)- 15 q y^2\Big)a_0^0,\nonumber \\
\frac{8}{105}e^{-\frac{q^2}{2}}&{}\Big(q^7 (x + i y)^3 - 21 i q^5 (x + i y)^2 y - 105 q^3 (x + i y) y^2 + 
 105 i q y^3\Big)a_0^0.
\end{align}
Another observation is that for a chosen homogeneity $m$ in $z, \bar{z}$, the highest exponent of $q$ 
is at least $2m+1$ and our solution realizes this minimum. The representative 
symplectic monogenics valued in ${\mathcal S}_+$ 
were already given for each homogeneity in Lemma \ref{allevenmonogenics}.
\end{remark}
In the following Theorem, we characterize the solution space for $T_s$ separately in the even case 
(including both even powers of $X_s$ acting on ${\mathcal S}_+$ and odd powers of $X_s$ acting on 
${\mathcal S}_-$) and the odd case (including both odd powers of $X_s$ acting on ${\mathcal S}_+$ 
and even powers of $X_s$ acting on ${\mathcal S}_-$.)

\begin{theorem}\label{TRecurTheorem}
\begin{enumerate}
\item
The symplectic spinor of the homogeneity $m\in\mN_0$ in the variables $z,\bar{z}$,
\begin{equation}\label{TRecurVyrazq}
s=e^{-\frac{q^2}{2}}q\big(A^m(q)z^m+ A^{m-1}(q)z^{m-1}\bar{z}+ \ldots + A^1(q)z\bar{z}^{m-1}+A^0(q)\bar{z}^m \big)
\end{equation}
with coefficients in the formal power series in $q$,
$$
A^r=a^r_0+a^r_2 q^2+a^r_4 q^4 + \ldots,\, a^r_{k}\in\mC,\, r=0,\ldots ,m,\, k\in 2\mN_0,
$$ 
is in the kernel of the symplectic twistor operator $T_s$ provided 
the coefficients $a^r_k$ satisfy the recurrence relations
\begin{eqnarray}
& & 0=m k a^m_k+k a^{m-1}_k-2a^{m-1}_{k-2}, \nonumber\\
& & 0=(m-1)k a^{m-1}_k+2 k a^{m-2}_k-4a^{m-2}_{k-2}, \nonumber\\
& & \,\,\, \dots \nonumber \\
& & 0=2 k a^{2}_k+(m-1) k a^{1}_k-2(m-1)a^{1}_{k-2},\\ \nonumber
& & 0= k a^1_k+m k a^{0}_k-2ma^{0}_{k-2},
\end{eqnarray}
equivalent to
\begin{equation}\label{TRecurZkracenoSude}
(m-p) k a_k^{m-p}+(p+1) k a_k^{m-1-p}-2(p+1)a_{k-2}^{m-1-p}=0,
\end{equation}
for all $p=0,1,\ldots,m-1$.
\item
The symplectic spinor of the homogeneity $m\in\mN_0$ in the variables $z,\bar{z}$,
\begin{equation}\label{TRecurVyraz}
s=e^{-\frac{q^2}{2}}\big(A^m(q)z^m+ A^{m-1}(q)z^{m-1}\bar{z}+ \ldots + A^1(q)z\bar{z}^{m-1}+A^0(q)\bar{z}^m \big)
\end{equation}
with coefficients in the formal power series in $q$,
$$
A^r=a^r_0+a^r_2 q^2+a^r_4 q^4 + \ldots,\, a^r_{k}\in\mC,\, r=0,\ldots ,m,\, k\in 2\mN_0,
$$ 
is in the kernel of the symplectic twistor operator $T_s$ provided 
the coefficients $a^r_k$ satisfy the recurrence relations
\begin{eqnarray}
& & 0=m(k-1)a^m_k+(k-1)a^{m-1}_k-2a^{m-1}_{k-2}, \nonumber\\
& & 0=(m-1)(k-1)a^{m-1}_k+2(k-1)a^{m-2}_k-4a^{m-2}_{k-2}, \nonumber\\
& & \,\,\, \dots \nonumber \\
& & 0=2(k-1)a^{2}_k+(m-1)(k-1)a^{1}_k-2(m-1)a^{1}_{k-2},\\ \nonumber
& & 0=(k-1)a^1_k+m(k-1)a^{0}_k-2ma^{0}_{k-2},
\end{eqnarray}
equivalent to
\begin{equation}\label{TRecurZkraceno}
(m-p)(k-1)a_k^{m-p}+(p+1)(k-1)a_k^{m-1-p}-2(p+1)a_{k-2}^{m-1-p}=0,
\end{equation}
for all $p=0,1,\ldots,m-1$.
\end{enumerate}
\end{theorem}

{\bf Proof:}

Concerning the first part, we have
\begin{align*}
T_s\Big(e^{-\frac{q^2}{2}}q\big(&{}A^m(q)z^m+ A^{m-1}(q)z^{m-1}\bar{z}+ \ldots + A^1(q)z\bar{z}^{m-1}+A^0(q)\bar{z}^m \big)
\Big)\nonumber\\
=e^{-\frac{q^2}{2}}q^2 \Big(&{} z^{m-1}\left(m[-\partial_q ]A^m(q)+[2q - \partial_q]A^{m-1}(q)\right)\\
+&{}z^{m-2}\bar{z}\left((m-1)[-\partial_q] A^{m-1}(q)+2[2q - \partial_q]A^{m-2}(q)\right)\\
&{}\ldots \\
+&{}\bar{z}^m\left([-\partial_q] A^{1}(q)+m[2q - \partial_q]A^{0}(q)\right)\Big)=0,
\end{align*}
where  
\begin{eqnarray}
& & [-\partial_q]A^{r}(q)={}-2 a_2^r q-4a_4^rq^3-6a_6^rq^5-\ldots ,\nonumber \\ \nonumber 
& & [2q - \partial_q]A^{r}(q)={}(2a_0^r-2a_2^r)q+(2a_2^r-4a_4^r)q^3+ \ldots ,
\end{eqnarray}
etc. Then the coefficients of $A^r(q)=a^r_0+a^r_2 q^2+a^r_4 q^4 + \ldots$, $r=0,\ldots, m$ satisfy 
the recurrence relations
\begin{equation}\nonumber
(m-p)k a_k^{m-p}+(p+1)k a_k^{m-1-p}-2(p+1)a_{k-2}^{m-1-p}=0,\, p=0,\ldots, m-1.
\end{equation}
As for the second part, we get
\begin{eqnarray}
& & (1-q\partial_q-q^2)e^{-\frac{q^2}{2}}A^r(q)=e^{-\frac{q^2}{2}}[1-q\partial_q]A^r(q),
\nonumber \\
& & (1-q\partial_q+q^2)e^{-\frac{q^2}{2}}A^r(q)=e^{-\frac{q^2}{2}}[1+2q^2-q\partial_q]A^r(q).
\end{eqnarray}
The annihilation condition for the symplectic twistor operator $T_s$ acting on (\ref{TRecurVyraz}) is equivalent to
\begin{align}
T_s\Big(e^{-\frac{q^2}{2}}\big(&{}A^m(q)z^m+ A^{m-1}(q)z^{m-1}\bar{z}+ \ldots + A^1(q)z\bar{z}^{m-1}+A^0(q)\bar{z}^m \big)
\Big)\nonumber\\
=e^{-\frac{q^2}{2}}\Big(&{}
z^{m-1}\big(m[1-q\partial_q]A^m(q)+[1+2q^2-q\partial_q]A^{m-1}(q)\big)\nonumber\\
&{}z^{m-2}\bar{z}\big((m-1)[1-q\partial_q]A^{m-1}(q)+2[1+2q^2-q\partial_q]A^{m-2}(q)\big)\nonumber\\
&{}\vdots\nonumber\\
&{}z\bar{z}^{m-1}\big(2[1-q\partial_q]A^2(q)+(m-1)[1+2q^2-q\partial_q]A^1(q)\big)\nonumber\\
&{}\bar{z}^m\big([1-q\partial_q]A^{1}(q)+m[1+2q^2-q\partial_q]A^0(q)\big)
\Big),
\end{align}
and this completes the proof of the Theorem.

\hfill $\square$

\begin{remark}
The explicit solution vectors for the symplectic twistor operator $T_s$ 
are, for the choice of $A^0(q)=a_0^0\neq 0$, given in homogeneities $m=1,2,3$ by  
\begin{align}
e^{-\frac{q^2}{2}}\big(&{}\left( -1 + 2 q^2\right) z + \bar{z}\big)a_0^0, 
\nonumber \\ \nonumber 
e^{-\frac{q^2}{2}}\bigg(&{}\Big(1 - 4 q^2 + \frac{4}{3} q^4\Big) z^2 + \left(-2 + 4 q^2\right) z \bar{z} + \bar{z}^2\bigg)a_0^0,
\\ \nonumber
e^{-\frac{q^2}{2}}\bigg(&{}\Big(-1 + 6 q^2 - 4 q^4 + \frac{8}{15} q^6\Big) z^3 + \left(3 - 12 q^2 + 4 q^4\right) z^2 \bar{z} 
\\ \nonumber
&{}+ \left(-3 + 6 q^2\right) z \bar{z}^2 + \bar{z}^3\bigg)a_0^0.\label{TPrikHom3}
\end{align}
The same solutions expressed in the variables $x,y$ are
\begin{align}
2 e^{-\frac{q^2}{2}}&{}\left( q^2 (x + i y) - i y\right)a_0^0,\nonumber \\
\frac{4}{3}e^{-\frac{q^2}{2}}&{}\left(q^4 (x + i y)^2 + 6 q^2 y (-i x + y)- 3 y^2\right)a_0^0,\nonumber \\
\frac{8}{15}e^{-\frac{q^2}{2}}&{}\left(q^6 (x + i y)^3 - 15 i q^4 (x + i y)^2 y - 45 q^2 (x + i y) y^2 + 
 15 i y^3\right)a_0^0.
\end{align}
\end{remark}

\begin{theorem}
Let $s=s(z,\bar z,q)\in \mathrm{Pol}(\mathbb{R}^2,{\mathcal S}_-)$ be a polynomial symplectic spinor in the solution space
of the symplectic Dirac operator $D_s$, i.e. the symplectic spinor $s$ satisfying the recurrence relations in the first part of 
Theorem (\ref{DsRecurTheorem}). Then $X_s(s)$ is in kernel of the symplectic twistor operator, $T_s\big(X_s(s)\big)=0$.
\end{theorem}
{\bf Proof:}

Let us consider the polynomial symplectic spinor of homogeneity $m$,
$$
s=e^{-\frac{q^2}{2}}q\big(A^m(q)z^m+ A^{m-1}(q)z^{m-1}\bar{z}+ \ldots + A^1(q)z\bar{z}^{m-1}+A^0(q)\bar{z}^m \big),
$$ 
where $A^r(q)=a^r_0+a^r_2 q^2+a^r_4 q^4 + \ldots$, $r=0,\ldots,m$ satisfies the recursive relations 
(\ref{DsRecurZkraceno}). The functions $A^l(q)qe^{-\frac{q^2}{2}}$, $l=0,\ldots,m$ are in the class 
of Schwartz functions. We use the notational simplification $s(z,\bar z,q)=e^{-\frac{q^2}{2}}qW$, $W=W(z,\bar z,q)$.
Then
$$
X_s(e^{-\frac{q^2}{2}}qW)=e^{-\frac{q^2}{2}}\left((2q^2-1-q\partial_q) z W +(1+q\partial_q)\bar{z}W\right),
$$
which can be rewritten as
$$
X_s(e^{-\frac{q^2}{2}}qW)=e^{-\frac{q^2}{2}}\big(B^{m+1}(q)z^{m+1}+ B^{m}(q)z^{m}\bar{z}+ \ldots +B^0(q)\bar{z}^{m+1} \big),
$$
where $B^r(q)=b^r_0+b^r_2 q^2+b^r_4 q^4 + \ldots$, $r=0,\ldots,m+1$, and the coefficients of this formal power series satisfy 
\begin{equation}\label{vztahba}
b_k^m=2a_{k-2}^{m-1}+(k+1)(a_k^m-a_k^{m-1}).
\end{equation}
We show that $B^r(q)$ satisfy the recurrence relations (\ref{TRecurZkraceno}) for $p=0,1,\ldots,m$ in Theorem (\ref{TRecurTheorem}).
It follows from (\ref{vztahba}) that
{\footnotesize
\begin{align}
&{} (m+1-p)(k-1)\big(2a_{k-2}^{m-p}+(k+1)(a_k^{m-p+1}-a_k^{m-p}) \big)
\nonumber \\
&{} +(p+1)(k-1)\big(2a_{k-2}^{m-p-1}+(k+1)(a_k^{m-p}-a_k^{m-p-1}) \big)
\nonumber \\
&{} -2(p+1)\big(2a_{k-4}^{m-p-1}+(k-1)(a_{k-2}^{m-p}-a_{k-2}^{m-p-1}) \big)
\nonumber \\
&{} =2\big((m-p)(k-1)a_{k-2}^{m-p}+(p+1)(k-1)a_{k-2}^{m-p-1}-2(p+1)a_{k-4}^{m-p-1} \big)
\nonumber \\
&{} +(k-1)\big((m-p+1)(k+1)a_k^{m-p+1}+p(k+1)a_k^{m-p}-2pa_{k-2}^{m-p} \big)
\nonumber \\
&{} -(k-1)\big((m-p)(k-1)a_{k}^{m-p}+(p+1)(k-1)a_{k}^{m-p-1}-2(p+1)a_{k-2}^{m-p-1} \big)
\nonumber \\
&{} +2(k-1)a_{k-2}^{m-p}-(k-1)(k+1)a_k^{m-p}+(k-1)(k+1)a_k^{m-p}-2(k-1)a_{k-2}^{m-p} \nonumber \\
&{}  =0,
\end{align}}
where we used for the last equality the relation (\ref{DsRecurZkraceno}) to verify 
that each of the three rows in the last but one expression equals to zero.
The proof is complete.

\hfill
$\square$

\begin{theorem}
\label{notinkert}
Let $s=s(z,\bar z,q)\in \mathrm{Pol}_{m}(\mathbb{R}^2,{\mathcal S}_-)$, $m\in \mN_0$, be a symplectic spinor polynomial in the solution space
of the symplectic Dirac operator $D_s$. Then $s$ is not in the kernel of the symplectic twistor operator $T_s$ if and only if $m>0$.
\end{theorem}
{\bf Proof:}

By our assumption, the symplectic spinor $s$ satisfies the recurrence relation in Theorem \ref{DsRecurTheorem}.
Recall the recurrence relations for symplectic spinors valued in $\mathcal{S}_-$, which are in the 
solution space of $\mathrm{Ker} (T_s)$, (\ref{TRecurZkracenoSude}):
\begin{equation}\nonumber
(m-p)k a_k^{m-p}+(p+1)k a_k^{m-1-p}-2(p+1)a_{k-2}^{m-1-p}=0,\, p=0,\ldots, m-1.
\end{equation}
By Theorem \ref{DsRecurTheorem}, the coefficients $a_k^{r}$ satisfy the relations (\ref{DsRecurZkraceno}) 
\begin{equation*}
(m-p)(k+1)a_k^{m-p}+(p+1)(k+1)a_k^{m-1-p}+2(p+1)a_{k-2}^{m-1-p}=0.
\end{equation*}
The comparison of the last two relations leads to 
\begin{equation}\label{Tphiproof}
(m-p)a_k^{m-p}+(p+1)a_k^{m-1-p}=0
\end{equation}
for all $k,p$, and these are just the coefficients by $q^{k+1}z^{m-1-p}\bar{z}^p$ in $T_s(s)$.
We choose the symplectic monogenic $s$ as in Remark \ref{remarkexplicitmonogenic}. 
For $k=2, p=0$, the coefficient in $T_s(s)$ by $q^3\bar{z}^{m-1}$ is  $\left( a_2^{1}+m a_2^{0}\right)$. 
Our choice for $s$ to be a solution for $D_s$ gives $a_2^1=\frac{2m}{3}a_0^0$ and $a_2^0=0$, 
therefore the coefficient in (\ref{Tphiproof}) will not be equal to zero and consequently will not be in 
$\mathrm{Ker} (T_s)$ for $m> 0$, $m\in\mN_0$. By 
$\mp(2,\mathbb{R})$-invariance, the whole metaplectic module does not belong to the kernel of $T_s$,
which finishes the proof.

\hfill 
$\square$

\begin{theorem}
\label{Ds2RecurRelace}
Let $m\in\mN_0, k\in2\mN_0$.   
\begin{enumerate}
\item
The recurrence relations for the coefficients $a^r_k$ of an even (even homogeneity in $q$) symplectic spinor $s$, 
\begin{equation*}
s=e^{-\frac{q^2}{2}}\big(A^m(q)z^m+ A^{m-1}(q)z^{m-1}\bar{z}+ \ldots + A^1(q)z\bar{z}^{m-1}+A^0(q)\bar{z}^m \big),
\end{equation*}
$A^r(q)=a^r_0+a^r_2 q^2+a^r_4 q^4 + \ldots$, $r=0,\ldots,m$, which is in the kernel of the square of the symplectic 
Dirac operator $D_s^2$, are 
\begin{align}\label{Ds2RecurRelZkracenoeven}
&{}(m-p)(m-p-1)(k+2)(k+1)a_{k+2}^{m-p}+\nonumber\\
&{}(m-1-p)(p+1)\big(2(k+2)(k+1)a_{k+2}^{m-1-p}-2(2k+1) a_k^{m-1-p}\big)+\nonumber\\
&{}(p+1)(p+2)\big((k+2)(k+1)a_{k+2}^{m-2-p}-2(2k+1)a_k^{m-2-p}+4a_{k-2}^{m-2-p}\big)\nonumber \\
&{}=0
\end{align}
for $p=0,\ldots,m-2$.
\item
The recurrence relations for the coefficients $a^r_k$ of an odd (odd homogeneity in $q$) symplectic spinor $s$, 
\begin{equation*}
s=e^{-\frac{q^2}{2}}q\big(A^m(q)z^m+ A^{m-1}(q)z^{m-1}\bar{z}+ \ldots + A^1(q)z\bar{z}^{m-1}+A^0(q)\bar{z}^m \big),
\end{equation*}
$A^r(q)=a^r_0+a^r_2 q^2+a^r_4 q^4 + \ldots$, $r=0,\ldots,m$, which is in the kernel of the square of the symplectic 
Dirac operator $D_s^2$, are 
\begin{align}\label{Ds2RecurRelZkracenoodd}
&{}(m-p)(m-p-1)(k+2)(k+3)a_{k+2}^{m-p}+\nonumber\\
&{}(m-1-p)(p+1)\big(2(k+2)(k+3)a_{k+2}^{m-1-p}-2(2k+3) a_k^{m-1-p}\big)+\nonumber\\
&{}(p+1)(p+2)\big((k+2)(k+3)a_{k+2}^{m-2-p}-2(2k+3)a_k^{m-2-p}+4a_{k-2}^{m-2-p}\big)\nonumber \\
&{}=0.
\end{align}
for $p=0,\ldots,m-2$.
\end{enumerate}
\end{theorem}
{\bf Proof:}

The second power of the symplectic Dirac operator $D_s$ is equal to 
\begin{equation}
D_s^2=(q^2+2q\partial_q+1+\partial_q^2)\partial_z^2+2(-q^2+\partial_q^2)\partial_z\partial_{\bar{z}}+(q^2-2q\partial_q -1+\partial_q^2)\partial_{\bar{z}}^2.
\end{equation}
In the even case, the action of $D_s^2$ results in
\begin{align}
D_s^2\Big( e^{-\frac{q^2}{2}}  \big( A^m&{}(q)z^m+ A^{m-1}(q)z^{m-1}\bar{z}+ \ldots +A^0(q)\bar{z}^m \big)\Big)
\nonumber \\
=e^{-\frac{q^2}{2}}\Big(z^{m-2}&{}\big( m(m-1)[\partial_q^2]A^m(q)+(m-1)[2\partial^2_q-4q\partial_q-2]A^{m-1}(q)
\nonumber \\
&{}+2[\partial_q^2-4q\partial_q-2+4q^2]A^{m-2} (q)\big){}+\ldots +
\nonumber \\
\bar{z}^{m-2}&{}\big( 2[\partial_q^2]A^2(q)+(m-1)[2\partial^2_q-4q\partial_q-2]A^{1}(q)+
\nonumber \\
&{}+ m(m-1)[\partial_q^2-4q\partial_q-2+4q^2]A^{0}(q) \big)
\Big),
\end{align}
where
\begin{align*}
[\partial_q^2]A^r(q)=\,& 2 a_2^r+ 12 a_4^r q^2 +\ldots \nonumber\\
[2\partial^2_q-4q\partial_q-2]A^{r}(q)=\,&4 a_2^r -2a_0^r + (24 a_4^r-8 a_2^r-2a_2^r)q^2+\ldots\nonumber\\
[\partial_q^2-4q\partial_q-2+4q^2]A^{r}(q)=\,&2a_2^r-2a_0^r+(12 a_4^r-8 a_2^r-2a_2^r+4a_0^r)q^2+\ldots
\end{align*}

The odd homogeneity case is analogous. Denoting $s=e^{-\frac{q^2}{2}}qW$, where $W=A^m(q)z^m+ \ldots +A^0(q)\bar{z}^m $, we get
\begin{align*}
\partial_z^2(q^2+2q\partial_q+1+\partial_q^2)e^{-\frac{q^2}{2}}qW=\,&{} \partial_z^2 e^{-\frac{q^2}{2}}[2\partial_q+q\partial_q^2]W,\\
2\partial_z\partial_{\bar{z}}(-q^2+\partial_q^2)e^{-\frac{q^2}{2}}qW=\,&{} 2\partial_z\partial_{\bar{z}}e^{-\frac{q^2}{2}}[q\partial_q^2-2q^2\partial_q+2\partial_q-3q]W,\\
\partial_{\bar{z}}^2(q^2-2q\partial_q -1+\partial_q^2)e^{-\frac{q^2}{2}}qW=\,
&{} \partial_{\bar{z}}^2e^{-\frac{q^2}{2}}[q\partial_q^2-4q^2\partial_q+2\partial_q+4q^3-6q]W,
\end{align*}
and the proof follows. 

The irreducible $\mp(2,\mR)$-submodules in the kernel of $D_s^2$ were 
put into boxes on the scheme of the $\mp(2,\mR)$-decomposition of $\mathrm{Pol}(\mR^2)\otimes{\mathcal S}$:

\begin{eqnarray}
\xymatrix@=11pt{
*+[F]{M_0} \ar[r] & *+[F]{X_s M_0} \ar @{} [d] |{\oplus} \ar[r] & X_s^2 M_0 \ar @{} [d] |{\oplus} \ar[r] & X_s^3 M_0 \ar @{} [d] |{\oplus}
 \ar[r] & X_s^4 M_0 \ar @{} [d] |{\oplus}\ar[r] & X_s^5 M_0 \ar @{} [d] |{\oplus} \\
& *+[F]{M_1} \ar[r] & *+[F]{X_s M_1} \ar @{} [d] |{\oplus}\ar[r] & X_s^2 M_1 \ar @{} [d] |{\oplus}
 \ar[r] & X_s^3 M_1 \ar @{} [d] |{\oplus}\ar[r] & X_s^4 M_1 \ar @{} [d] |{\oplus} \\
&& *+[F]{M_2} \ar[r] & *+[F]{X_s M_2} \ar @{} [d] |{\oplus}
 \ar[r] & X_s^2 M_2 \ar @{} [d] |{\oplus}\ar[r] & X_s^3 M_2 \ar @{} [d] |{\oplus} \\
&&& *+[F]{M_3} \ar[r] & *+[F]{X_s M_3} \ar @{} [d] |{\oplus}\ar[r] & X_s^2 M_3  \ar @{} [d] |{\oplus} \\
&&&& *+[F]{M_4} \ar[r] & *+[F]{X_s M_4} \ar @{} [d] |{\oplus}  \\
&&&&& *+[F]{M_5} 
}
\end{eqnarray}

\hfill
$\square$

\begin{theorem}
\label{inthekerneld2}
The solution space of the symplectic twistor operator $T_s$ is a subspace of the space of solutions of the square of 
the symplectic Dirac operator $D_s^2$. In particular,  
the recurrence relations for $D_s^2$ specialized to even resp. odd symplectic spinors from Theorem \ref{Ds2RecurRelace} 
are solved by (\ref{TRecurZkraceno}) resp. (\ref{TRecurZkracenoSude}).
\end{theorem}
{\bf Proof:}

Let us start with even symplectic spinors.
It is straigtforward to rewrite the recurrence relations in Theorem \ref{Ds2RecurRelace},
\begin{align*}
&{}(m-p)(m-p-1)(k+2)(k+1)a_{k+2}^{m-p}+\nonumber \\ \nonumber
&{}(m-1-p)(p+1)\big(2(k+2)(k+1)a_{k+2}^{m-1-p}-2(2k+1) a_k^{m-1-p}\big)+\\
&{}(p+1)(p+2)\big((k+2)(k+1)a_{k+2}^{m-2-p}-2(2k+1)a_k^{m-2-p}+4a_{k-2}^{m-2-p}\big)=0,\nonumber
\end{align*}
into
{\footnotesize
\begin{align*}
&{}(m-1-p)(k+2)\big((m-p)(k+1)a_{k+2}^{m-p}+(p+1)(k+1)a_{k+2}^{m-1-p}-2(p+1)a_{k}^{m-1-p} \big)\\
+&{}(p+1)(k+2)\big((m-1-p)(k+1)a_{k+2}^{m-1-p}+(p+2)(k+1)a_{k+2}^{m-2-p}-2(p+2)a_{k}^{m-2-p} \big)\\
-&{}2(p+1)\big((m-1-p)(k-1)a_{k}^{m-1-p}+(p+2)(k-1)a_{k}^{m-2-p}-2(p+2)a_{k-2}^{m-2-p} \big)=0.
\end{align*}}
Because each of the last three rows corresponds to a recurrence relation (\ref{TRecurZkraceno}), 
the claim follows.

In the odd case, the recurrence relations
\begin{align*}
&{}(m-p)(m-p-1)(k+2)(k+3)a_{k+2}^{m-p}+\nonumber \\ \nonumber
&{}(m-1-p)(p+1)\big(2(k+2)(k+3)a_{k+2}^{m-1-p}-2(2k+3) a_k^{m-1-p}\big)+\\
&{}(p+1)(p+2)\big((k+2)(k+3)a_{k+2}^{m-2-p}-2(2k+3)a_k^{m-2-p}+4a_{k-2}^{m-2-p}\big)=0,\nonumber
\end{align*}
can be rewritten as
{\footnotesize
\begin{align*}
&{}(m-1-p)(k+3)\big((m-p)(k+2)a_{k+2}^{m-p}+(p+1)(k+2)a_{k+2}^{m-1-p}-2(p+1)a_{k}^{m-1-p} \big)\\
+&{}(p+1)(k+3)\big((m-1-p)(k+2)a_{k+2}^{m-1-p}+(p+2)(k+2)a_{k+2}^{m-2-p}-2(p+2)a_{k}^{m-2-p} \big)\\
-&{}2(p+1)\big((m-1-p)k a_{k}^{m-1-p}+(p+2)k a_{k}^{m-2-p}-2(p+2)a_{k-2}^{m-2-p} \big)=0,
\end{align*}}
and each of the last three rows corresponds to the recurrence relation (\ref{TRecurZkracenoSude}).

\hfill
$\square$


\begin{theorem}
The solution space of the symplectic twistor operator $T_s$, acting on $\mathrm{Pol}(\mathbb{R}^2,{\mathcal S})$, consists
of the set of $\mp(2,\mR)$-modules pictured in the squares realized in the decomposition of
$\mathrm{Pol}(\mathbb{R}^2,{\mathcal S})$ on $\mp(2,\mR)$ irreducible subspaces, (\ref{obrdekonposition}):
\begin{enumerate}
\item
$\mathrm{Pol}(\mathbb{R}^2,{\mathcal S}_-)$:
\begin{eqnarray}
\xymatrix@=11pt{
*+[F]{M^-_0} \ar @{}[d] |{q e^{-\frac{q^2}{2}}} \ar[rr] && *+[F]{X_s M^-_0} \ar @{} [d] |{\oplus} \ar[rr] && X_s^2 M^-_0 \ar @{} [d] |{\oplus} \ar[rr] && X_s^3 M^-_0 \ar @{} [d] |{\oplus}
  \ar[rr] &&  \ldots\\
&& M^-_1\ar @{}[d] |{e^{-\frac{q^2}{2}}(q^3 (x + i y) - 3 i q y)} \ar[rr] && *+[F]{X_s M^-_1} \ar @{} [d] |{\oplus}\ar[rr] && X_s^2 M^-_1 \ar @{} [d] |{\oplus}
  \ar[rr] && \ldots\\
&&&\ar @{}[d] |{ e^{-\frac{q^2}{2}}( q^5 (x + i y)^2 + 10 q^3 y (-i x + y) - 15 q y^2)}& M^-_2  \ar[rr] && *+[F]{X_s M^-_2} \ar @{} [d] |{\oplus}
  \ar[rr] && \ldots\\
&&&&&& M^-_3 \ar[rr] &&\ldots \\
}
\end{eqnarray}

\item
$\mathrm{Pol}(\mathbb{R}^2,{\mathcal S}_+)$:
\begin{eqnarray}
\xymatrix@=11pt{
*+[F]{M^+_0} \ar @{}[d] |{e^{-\frac{q^2}{2}}} \ar[rr] && *+[F]{X_s M^+_0} \ar @{} [d] |{\oplus} \ar[rr] && X_s^2 M^+_0 \ar @{} [d] |{\oplus} \ar[rr] && X_s^3 M^+_0 \ar @{} [d] |{\oplus}
  \ar[rr] &&  \ldots\\
&& M^+_1\ar @{}[d] |{e^{-\frac{q^2}{2}}(x + i y)} \ar[rr] && *+[F]{X_s M^+_1} \ar @{} [d] |{\oplus}\ar[rr] && X_s^2 M^+_1 \ar @{} [d] |{\oplus}
  \ar[rr] && \ldots\\
&&&&\ar @{}[d] |{ e^{-\frac{q^2}{2}} (x + i y)^2 } M^+_2  \ar[rr] && *+[F]{X_s M^+_2} \ar @{} [d] |{\oplus}
  \ar[rr] && \ldots\\
&&&&&& M^+_3 \ar[rr] &&\ldots \\
}
\end{eqnarray}

\end{enumerate}
Notice that the representative vectors in the solution space of $D_s$ are pictured under the spaces
of symplectic monogenics. In the case of ${\mathcal S}_+$, we exploit the symplectic 
monogenics constructed in Theorem \ref{stevenelement}.
\end{theorem}
{\bf Proof:}

It follows from the metaplectic Howe duality, \cite{bss}, that Theorem \ref{inthekerneld2} 
characterizes the $\mp(2,\mR)$-submodule of $\mathrm{Pol}(\mathbb{R}^2,{\mathcal S})$ contained in the solution space of $T_s$.
Then Theorem \ref{stevenelement}, Theorem \ref{notinkert} and Theorem \ref{Ds2RecurRelace} characterize the space 
of solutions as the image of the space of symplectic monogenics by $X_s$, in addition to the space 
of constant symplectic spinors. The proof is complete.

\hfill
$\square$

In previous sections, we discussed the space of polynomial solutions. A natural
question is an extension of the function space from polynomials to the class of 
analytic, smooth, Frechet, hyperfunction, generalized, etc.,
function spaces. For example, one can consider convergent power series constructed from the polynomial solutions. We 
shall not attempt to discuss this question in a greater generality, but observe 
the existence of solutions in wider classes of function spaces. 
 
Let us consider the function $z^n f(q)$ for $f\in S(\mR)$, $n\in\mN_0$. The substitution into (\ref{rce}) implies 
that it belongs to the solution space of $T_s$ provided  
$f(q)$ solves the ordinary differential equation
\begin{equation}\label{holrov}
(1-q^2)f(q)=q\frac{\partial}{\partial q} f(q).
\end{equation}
This equation has a unique solution $f(q)=q e^{- \frac{q^2}{2}}$ in $S(\mR)$, and so
$z^n q e^{- \frac{q^2}{2}}$ 
are in the kernel of the symplectic twistor operator for all $n\in \mathbb{N}_0$.

A generalization of this result is contained in the following lemma. 

\begin{lemma}
Let $h(z)$ be arbitrary holomorphic function on $\mC$.
Then the complex analytic symplectic spinor 
\begin{eqnarray}
h(z)q e^{-\frac{q^2}{2}}
\end{eqnarray} 
is in the kernel of the symplectic twistor operator $T_s$.
\end{lemma}
Consequently, the space of holomorphic functions on $\mC$ is embedded
into the space of smooth solutions of the symplectic twistor operator $T_s$.

Notice that an admissible continuous representation spaces of a reductive Lie group $G$ can 
be conveniently described in terms of a globalization of the underlying Harish-Chandra 
$(\gog,K)$-module, where $\gog$ resp. $K$ are the Lie algebra resp. the maximal compact 
subgroup of $G$. In this way, one has continuous representations of $G$ on the space 
of analytic, smooth, Frechet, hyperfunction, generalized, etc., functions. However, in our case of
$G=\Mp(2,\mR)$, $\gog=\mp(2,\mR)$ and $\tilde{K}$ given by the twofold covering of $\mathrm{U}(1)$, the representation 
on symplectic spinors is not admissible - in its composition series there are infinite
multiplicities of certain $G$-representations.
This means that the functional analytic tools developed in the representation theory are not straightforward to
apply here. On the other hand, it is still natural to ask for a characterization of the space of 
analytic, smooth, Frechet, hyperfunction, generalized, etc., solutions of both $T_s$ and $D_s$.

\vspace{0.5cm}
{\bf Acknowledgement:}  The authors gratefully acknowledge the support of the grant GA CR P201/12/G028 and SVV-2012-265317.

\vspace{0.3cm}

Marie Dost\'alov\'a, Petr Somberg

Mathematical Institute of Charles University,

Sokolovsk\'a 83, Praha 8 - Karl\'{\i}n, Czech Republic, 

E-mail: madost@seznam.cz, somberg@karlin.mff.cuni.cz.

\end{document}